\documentclass[12pt,vatola]{article}

\usepackage{graphicx}

\def\c{\centerline}

\def\re#1{\par\hangindent\parindent\indent\llap{#1\enspace}\ignorespaces}

\def\no{\noindent}

\begin{document}

\no{{\bf e-print:} {\it arXiv: math.GM/0606702}}

\vskip 15mm

\c{\large\bf Combinatorial Speculations and}\vskip 3mm

\c{\large\bf the Combinatorial Conjecture for Mathematics} \vskip
5mm

\c{Linfan Mao} \vskip 3mm \c{\scriptsize (Chinese Academy of
Mathematics and System Sciences, Beijing 100080)}

\c{\scriptsize maolinfan@163.com}

\vskip 8mm
\begin{minipage}{130mm}
\no{\bf Abstract}: {\small Combinatorics is a powerful tool for
dealing with relations among objectives mushroomed in the past
century. However, an more important work for mathematician is to
apply combinatorics to other mathematics and other sciences not
merely to find combinatorial behavior for objectives. Recently,
such research works appeared on journals for mathematics and
theoretical physics on cosmos. The main purpose of this paper is
to survey these thinking and ideas for mathematics and
cosmological physics, such as those of multi-spaces, map
geometries and combinatorial cosmoses, also the combinatorial
conjecture for mathematics proposed by myself in 2005. Some open
problems are included for the 21th mathematics by a combinatorial
speculation.}

\vskip 2mm \no{\bf Key words:} {\small  combinatorial speculation,
combinatorial conjecture for mathematics, Smarandache multi-space,
M-theory, combinatorial cosmos.}

 \vskip 2mm \no{{\bf
Classification:} AMS(2000) 03C05,05C15,51D20,51H20,51P05,83C05,
83E50.}
\end{minipage}

\vskip 8mm

\no{\bf $1$. The role of classical combinatorics in mathematics}

\vskip 6mm

\no Modern science has so advanced that to find a universal genus
in the society of sciences is nearly impossible. Thereby a
scientist can only give his or her contribution in one or several
fields. The same thing also happens for researchers in
combinatorics. Generally, combinatorics deals with twofold
questions:

\vskip 3mm

{\bf Question $1.1$.}\ {\it determine or find structures or
properties of configurations, such as those structure results
appeared in graph theory, combinatorial maps and design
theory,..., etc..}

\vskip 2mm

{\bf Question $1.2$.}\ {\it enumerate configurations, such as
those appeared in the enumeration of graphs, labelled graphs,
rooted maps, unrooted maps and combinatorial designs,...,etc..}

\vskip 2mm

Consider the contribution of a question to science. We can
separate mathematical questions into three grades:

\vskip 3mm

{\bf Grade $1$} \ {\it they contribute to all sciences.}

\vskip 2mm

{\bf Grade $2$} \ {\it they contribute to all or several branches
of mathematics.}

\vskip 2mm

{\bf Grade $3$} \ {\it they contribute only to one branch of
mathematics, for instance, just to the graph theory or
combinatorial theory.}

\vskip 2mm

Classical combinatorics is just a {\it grade $3$ mathematics} by
this view. This conclusion is gloomy for researchers in
combinatorics, also for me 4 years ago. {\it Whether can
combinatorics be applied to other mathematics or other sciences?
Whether can it contribute to human's lives, not just in games?}

Although become a universal genus in science is nearly impossible,
{\it our world is a combinatorial world}. A combinatorician should
stand on all mathematics and all sciences, not just on classical
combinatorics and then with a real combinatorial notion, i.e.,
{\it combining different fields into a unifying field}
([25]-[28]), such as combine different or even anti branches in
mathematics or science into a unifying science for its freedom of
research ([24]). This combination also requires us answering three
questions for solving a combinatorial question before. {\it What
is this question working for? What is its objective? What is its
contribution to science or human's society?} After these works be
well done, modern combinatorics can be applied to all sciences and
all sciences are combinatorization.

\vskip 5mm

\no{\bf $2$. The combinatorics metrization and mathematics
combinatorization}

\vskip 4mm

\no There is a prerequisite for the application of combinatorics
to other branch mathematics and other sciences, i.e, to introduce
various metrics into combinatorics, ignored by the classical
combinatorics since they are the fundamental of scientific
realization for our world. This speculation is firstly appeared in
the beginning of Chapter $5$ of my book $[16]$:

\vskip 3mm

{\it  $\cdots$ our world is full of measures. For applying
combinatorics to other branch of mathematics, a good idea is
pullback measures on combinatorial objects again, ignored by the
classical combinatorics and reconstructed or make combinatorial
generalization for the classical mathematics, such as those of
algebra, differential geometry, Riemann geometry, Smarandache
geometries, $\cdots$ and the mechanics, theoretical physics,
$\cdots$. }\vskip 2mm

The combinatorial conjecture for mathematics, abbreviated to {\it
CCM} is stated in the following.

\vskip 4mm

\no{\bf Conjecture $2.1$}(CCM Conjecture) \ {\it Mathematics can
be reconstructed from or turned into combinatorization.}

\vskip 3mm

\no{\bf Remark $2.1$} \ We need some further clarifications for
this conjecture.

\vskip 3mm

($i$)\ \ \ This conjecture assumes that one can selects finite
combinatorial rulers and axioms to reconstruct or make
generalization for classical mathematics.

($ii$) \ Classical mathematics is a particular case in the
combinatorization of mathematics, i.e., the later is a
combinatorial generalization of the former.

($iii$)\ We can make combinatorizations of different branches in
mathematics into one and find new theorems after then.

\vskip 2mm

Therefore, a branch in mathematics can not be ended if it has not
been combinatorization and all mathematics can not be ended if its
combinatorization has not completed. There is an assumption in
one's realization of our world, i.e., {\it every science can be
turned into mathematization}. Whence, we similarly get the
combinatorial conjecture for sciences.

\vskip 4mm

\no{\bf Conjecture $2.2$}(CCS Conjecture) \ {\it Sciences can be
reconstructed from or turned into combinatorization.}

\vskip 3mm

A typical example for the combinatorization of classical
mathematics is the {\it combinatorial map theory}, i.e., a
combinatorial theory for surfaces([14]-[15]). Combinatorially, a
surface is topological equivalent to a polygon with even number of
edges by identifying each pairs of edges along a given direction
on it. If label each pair of edges by a letter $e,e\in {\mathcal
E}$, a surface $S$ is also identifying with a cyclic permutation
such that each edge $e,e\in {\mathcal E}$ just appears two times
in $S$, one is $e$ and another is $e^{-1}$. Let $a,b,c,\cdots$
denote the letters in ${\mathcal E}$ and $A,B,C,\cdots$ the
sections of successive letters in a linear order on a surface $S$
(or a string of letters on $S$). Then, a surface can be
represented as follows:

$$S=(\cdots , A,a,B,a^{-1},C,\cdots),$$

\no{where, $a\in {\mathcal E}$,$A,B,C$ denote a string of letters.
Define three elementary transformations as follows:}

\vskip 2mm $(O_1)\quad\quad (A,a,a^{-1},B)\Leftrightarrow (A,B);$

\vskip 2mm $(O_2)\quad\quad (i)\quad
(A,a,b,B,b^{-1},a^{-1})\Leftrightarrow (A,c,B,c^{-1}) ;$

\quad\quad\quad\quad $(ii)\quad (A,a,b,B,a,b)\Leftrightarrow
(A,c,B,c); $

\vskip 2mm $(O_3)\quad\quad (i)\quad
(A,a,B,C,a^{-1},D)\Leftrightarrow (B,a,A,D,a^{-1},C);$

\quad\quad\quad\quad $(ii)\quad (A,a,B,C,a,D)\Leftrightarrow
(B,a,A,C^{-1},a,D^{-1}).$ \vskip 2mm

If a surface $S$ can be obtained from $S_0$ by these elementary
transformations $O_1$-$O_3$, we say that $S$ is elementary
equivalent with $S_0$, denoted by $S\sim_{El}S_0$. Then we can get
the classification theorem of compact surfaces as follows($[29]$):

\vskip 3mm

{\it Any compact surface is homeomorphic to one of the following
standard surfaces:}

($P_0$) {\it the sphere: $aa^{-1}$;}

($P_n$) {\it the connected sum of $n,n\geq 1$ tori:}

$$a_1b_1a_1^{-1}b_1^{-1}a_2b_2a_2^{-1}b_2^{-1}\cdots a_nb_na_n^{-1}b_n^{-1};$$

($Q_n$) {\it the connected sum of $n,n\geq 1$ projective planes:}

$$a_1a_1a_2a_2\cdots a_na_n.$$

\vskip 2mm

A {\it map} $M$ is a connected topological graph cellularly
embedded in a surface $S$. In 1973, Tutte suggested an algebraic
representation for an embedding graph on a locally orientable
surface $([16])$:

A {\it{combinatorial map}} $M = ({\cal X}
_{\alpha,\beta},\cal{P})$ is defined to be a basic permutation
$\cal{P}$, i.e, for any $x\in {\cal X}_{\alpha,\beta}$, no integer
$k$ exists such that ${\cal{P}}^{k}x = \alpha x$, acting on ${\cal
X} _{\alpha,\beta}$, the disjoint union of {\it quadricells} $Kx$
of $x\in  X$ (the base set), where $K=\{1,\alpha,\beta,\alpha\beta
\}$ is the {\it Klein group} satisfying the following two
conditions:\vskip 3mm

($i$) \ \ {\it  $\alpha{\cal{P}}={\cal{P}}^{-1}\alpha$;}\vskip 2mm

($ii$) \ {\it the group $\Psi_{J}=<\alpha,\beta,\cal{P}>$ is
transitive on ${\cal X}_{\alpha,\beta}$.}\vskip 2mm

For a given map $M=({\mathcal X}_{\alpha ,\beta},{\mathcal P})$,
it can be shown that $M^* = ({\mathcal X}_{\beta ,\alpha
},{\mathcal P}\alpha\beta)$ is also a map, call it the {\it dual}
of the map $M$. The vertices of $M$ are defined as the pairs of
conjugatcy orbits of ${\mathcal P}$ action on ${\mathcal
X}_{\alpha,\beta}$ by the condition $(i)$ and edges the orbits of
$K$ on ${\mathcal X}_{\alpha,\beta}$, for example, for $\forall
x\in {\mathcal X}_{\alpha,\beta}$, $\{x,\alpha x,\beta
x,\alpha\beta x\}$ is an edge of the map $M$. Define the faces of
$M$ to be the vertices in the dual map $M^*$. Then the Euler
characteristic $\chi (M)$ of the map $M$ is

$$\chi (M)= \nu (M)-\varepsilon (M)+\phi (M)$$

\no{where, $\nu (M), \varepsilon (M), \phi (M)$ are the number of
vertices, edges and faces of the map $M$, respectively. For each
vertex of a map $M$, its valency is defined to be the length of
the orbits of ${\mathcal P}$ action on a quadricell incident with
$u$. }

For example, the graph $K_4$ on the tours with one face length $4$
and another $8$ shown in Fig.$2.1$

\vskip 2mm

\includegraphics[bb=-10 10 200 200]{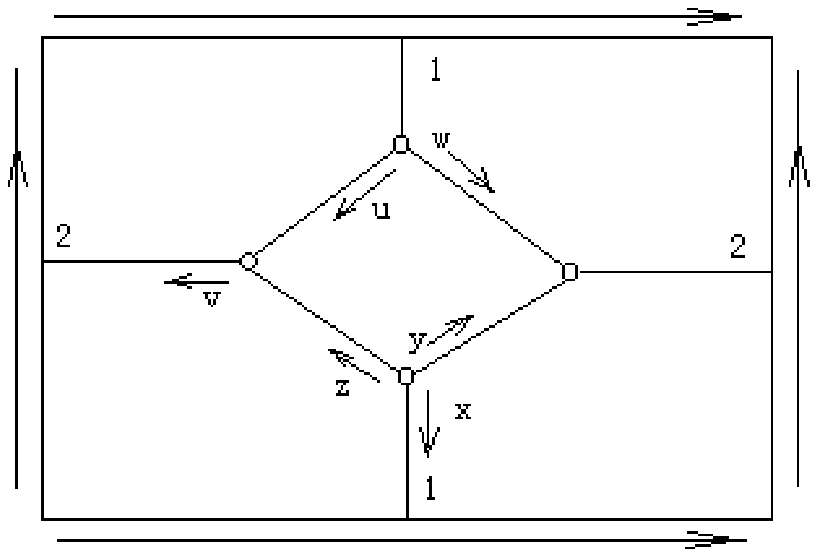}

\vskip 2mm \c{\bf Fig.$2.1$} \vskip 3mm

\no can be algebraically represented by $({\mathcal
X}_{\alpha,\beta},\mathcal{P})$ with ${\mathcal X}_{\alpha,\beta}=
\{x,y,z,u,v,w,\alpha x,\alpha y,\alpha z,$ $ \alpha u,\alpha
v,\alpha w, \beta x,\beta y,\beta z,\beta u,
 \beta v,\beta w,\alpha\beta x,
\alpha \beta y,\alpha \beta z,\alpha \beta u,\alpha \beta
v,\alpha\beta w \}$ and

\begin{eqnarray*}
{\mathcal P} &=& (x,y,z)(\alpha \beta x,u,w)(\alpha \beta z,\alpha
\beta u,v)
(\alpha \beta y,\alpha \beta v,\alpha \beta w)\\
&\times& (\alpha x,\alpha z,\alpha y)(\beta x,\alpha w,\alpha
u)(\beta z,\alpha v,\beta u)(\beta y,\beta w,\beta v)
\end{eqnarray*}

\no with $4$ vertices, $6$ edges and $2$ faces on an orientable
surface of genus $1$.

By the view of combinatorial maps, these standard surfaces $P_0,
P_n, Q_n$ for $n\geq 1$ is nothing but the bouquet $B_n$ on a
locally orientable surface with just one face. Therefore,
combinatorial maps are the combinatorization of surfaces.

Many open problems are motivated by the {\it CCM Conjecture}. For
example,  a {\it Gauss mapping} among surfaces is defined as
follows.\vskip 3mm

{\it Let ${\mathcal S}\subset R^3$ be a surface with an
orientation {$\bf N$}. The mapping {$\bf N:$}${\mathcal
S}\rightarrow R^3$ takes its value in the unit sphere}

$$S^2=\{(x,y,z)\in R^3|x^2+y^2+z^2=1\}$$

\no{\it along the orientation {$\bf N$}. The map ${\bf N:}
{\mathcal S}\rightarrow S^2$, thus defined, is called the Gauss
mapping.}\vskip 3mm

we know that for a point $P\in {\mathcal S}$ such that the
Gaussian curvature $K(P)\not=0$ and $V$ a connected neighborhood
of $P$ with $K$ does not change sign,

$$K(P)=\lim_{A\rightarrow 0}\frac{N(A)}{A},$$

\no{where $A$ is the area of a region $B\subset V$ and $ N(A)$ is
the area of the image of $B$ by the Gauss mapping $ N: {\mathcal
S}\rightarrow S^2$([2],[4]). Questions for the Gauss mapping
are}\vskip 3mm

$(i)$ {\it what is its combinatorial meaning of the Gauss mapping?
How to realizes it by combinatorial maps?}

$(ii)$ {\it how can we define various curvatures for maps and
rebuilt the results in the classical differential geometry?}\vskip
2mm

Let ${\mathcal S}$ be a compact orientable surface. Then the {\it
Gauss-Bonnet theorem} asserts that

$$\int\int_{\mathcal S}Kd\sigma =2\pi\chi ({\mathcal S}),$$

\no{where $K$ is the Gaussian curvature of ${\mathcal S}$.}

By the {\it CCM Conjecture}, the following questions should be
considered.\vskip 3mm

$(i)$  {\it How can we define various metrics for combinatorial
maps, such as those of length, distance, angle, area,
curvature,$\cdots$? }

($ii$) {\it Can we rebuilt the Gauss-Bonnet theorem by maps for
dimensional $2$ or higher dimensional compact manifolds without
boundary?}\vskip 2mm

One can see the references $[15]$ and $[16]$ for more open
problems for the classical mathematics motivated by this {\it CCM
Conjecture}, also raise new open problems for his or her research
works.

\vskip 5mm

\no{\bf $3$. The contribution of combinatorial speculation to
mathematics}

\vskip 4mm

\no{\bf $3.1$. The combinatorization of algebra}

\vskip 3mm

\no By the view of combinatorics, algebra can be also seen as a
combinatorial mathematics. The combinatorial speculation can
generalize it by means of the combinatorization. For this
objective, these Smarandache multi-algebraic systems are needed,
defined in the following.

\vskip 4mm

\no{\bf Definition $3.1$}([17],[18])\ {\it For any integers $n,
n\geq 1$ and $i, 1\leq i \leq n$, let $A_i$ be a set with an
operation set $O(A_i)$ such that $(A_i, O(A_i))$ is a complete
algebraic system. Then the union}

$$\bigcup\limits_{i=1}^n(A_i, O(A_i))$$

\no{\it is called an $n$ multi-algebra system.}

\vskip 3mm

An example of multi-algebra systems can be constructed by a finite
additive group. Let $n$ be an integer,
$Z_1=(\{0,1,2,\cdots,n-1\},+)$ an additive group $(mod {\rm n})$
and $P=(0,1,2,\cdots,n-1)$ a permutation. For any integer $i,
0\leq i\leq n-1$, define

$$Z_{i+1}=P^i(Z_1)$$

\no such that $P^i(k)+_iP^i(l)=P^i(m)$ in $Z_{i+1}$ if $k+l=m$ in
$Z_1$, where $+_i$ denotes the binary operation
$+_i:(P^i(k),P^i(l))\rightarrow P^i(m)$. Then we know that

$$\bigcup\limits_{i=1}^nZ_i$$

\no is an $n$ multi-algebra system .

The conception of multi-algebra systems can be extensively applied
for generalizing conceptions and results in the algebraic
structure, such as those of groups, rings, bodies, fields and
vector spaces, $\cdots$, etc.. Some of them are explained in the
following.

\vskip 4mm

\no{\bf Definition $3.2$} \ {\it Let
$\widetilde{G}=\bigcup\limits_{i=1}^n G_i$ be a complete
multi-algebra system with a binary operation set
$O(\widetilde{G})=\{\times_i, 1\leq i\leq n\}$. If for any integer
$i, 1\leq i\leq n$, $(G_i;\times_i)$ is a group and for $\forall
x,y,z\in \widetilde{G}$ and any two binary operations ¡°$\times$¡±
and ¡°$\circ$¡±, $\times\not= \circ$, there is one operation, for
example the operation $\times$ satisfying the distribution law to
the operation ¡°$\circ$¡± provided their operation results exist ,
i.e.,

$$x\times (y\circ z) = (x\times y)\circ (x\times z),$$

$$(y\circ z)\times x = (y\times x)\circ (z\times x),$$

\no then $\widetilde{G}$ is called a multi-group.}

\vskip 3mm

For a multi-group $(\widetilde{G},O(\widetilde{G}))$,
$\widetilde{G_1}\subset\widetilde{G}$ and
$O(\widetilde{G_1})\subset O(\widetilde{G})$, call
$(\widetilde{G_1},O(\widetilde{G_1}))$ a {\it sub-multi-group} of
$(\widetilde{G},O(\widetilde{G}))$ if $\widetilde{G_1}$ is also a
multi-group under the operations in $O(\widetilde{G_1})$, denoted
by $\widetilde{G_1} \preceq \widetilde{G}$. For two sets $A$ and
$B$, denote the union $A\bigcup B$ by $A\bigoplus B$ if $A\bigcap
B=\emptyset$. Then we get a generalization of the {\it Lagrange
theorem} of finite groups.

\vskip 4mm

\no{\bf Theorem $3.1$}([18]) \ {\it For any sub-multi-group
$\widetilde{H}$ of a finite multi-group $\widetilde{G}$, there is
a representation set $T$, $T\subset\widetilde{G}$, such that }

$$\widetilde{G}=\bigoplus\limits_{x\in T}x\widetilde{H}.$$

\vskip 3mm

For a sub-multi-group $\widetilde{H}$ of $\widetilde{G}$,
$\times\in O(\widetilde{H})$ and $\forall
g\in\widetilde{G}(\times)$, if for $\forall h\in\widetilde{H}$,

$$g\times h\times g^{-1}\in\widetilde{H},$$

\no then $\widetilde{H}$ is called a {\it normal sub-multi-group
of $\widetilde{G}$}. We call an arrangement of all operations in
$O(\widetilde{G})$ in order an {\it oriented operation sequence},
denote it by $\overrightarrow{O}(\widetilde{G})$. Then a
generalization of the {\it Jordan-H\"{o}lder theorem} for finite
multi-groups is described in the next result.

\vskip 4mm

\no{\bf Theorem $3.2$}([18]) \ {\it For a finite multi-group
$\widetilde{G}=\bigcup\limits_{i=1}^n G_i$ and an oriented
operation sequence $\overrightarrow{O}(\widetilde{G})$, the length
of maximal series of normal sub-multi-groups is a constant, only
dependent on $\widetilde{G}$ itself.}

\vskip 3mm

In Definition $2.2$, choose $n=2, G_1=G_2=\widetilde{G}$. Then
$\widetilde{G}$ is a body. If $(G_1;\times_1)$ and
$(G_2;\times_2)$ both are commutative groups, then $\widetilde{G}$
is a field. For multi-algebra systems with two or more operations
on one set, we introduce the conception of multi-rings and
multi-vector spaces in the following.

\vskip 4mm

\no{\bf Definition $3.3$} \ {\it Let
$\widetilde{R}=\bigcup\limits_{i=1}^mR_i$ be a complete
multi-algebra system with a double binary operation set
$O(\widetilde{R})=\{(+_i,\times_i) , 1\leq i\leq m\}$. If for any
integers $i, j, \ i\not= j, 1\leq i, j\leq m$, $(R_i; +_i,
\times_i)$ is a ring and for $\forall x,y,z\in\widetilde{R}$,

$$ (x+_iy)+_jz = x+_i(y+_jz), \  \ \ (x\times_iy)\times_jz = x\times_i(y\times_jz)$$

\no and

$$x\times_i(y+_jz) = x\times_iy +_jx\times_iz, \  \ \ (y+_jz)\times_ix = y\times_ix +_jz\times_ix$$

\no provided all thpse operation results exist, then
$\widetilde{R}$ is called a multi-ring. If for any integer $1\leq
i\leq m$, $(R;+_i,\times_i)$ is a filed, then $\widetilde{R}$ is
called a multi-filed.}

\vskip 4mm

\no{\bf Definition $3.4$} \ {\it Let
$\widetilde{V}=\bigcup\limits_{i=1}^k V_i$ be a complete
multi-algebra system with a binary operation set
$O(\widetilde{V})=\{ (\dot{+}_i,\cdot_i) \ | \ 1\leq i\leq m\}$
and $\widetilde{F}=\bigcup\limits_{i=1}^k F_i$ a multi-filed with
a double binary operation set $O(\widetilde{F})=\{(+_i,\times_i )\
| \ 1\leq i\leq k\}$. If for any integers $i,j, \ 1\leq i, j\leq
k$ and $\forall {\bf a,b,c}\in\widetilde{V}$,
$k_1,k_2\in\widetilde{F}$,

$(i)$ $(V_i;\dot{+}_i,\cdot_i)$ is a vector space on $F_i$ with
vector additive $\dot{+}_i$ and scalar multiplication $\cdot_i$;

$(ii)$ $({\bf a}\dot{+}_i{\bf b})\dot{+}_j{\bf c}= {\bf
a}\dot{+}_i({\bf b}\dot{+}_j{\bf c})$;

$(iii)$ $(k_1+_i k_2)\cdot_j{\bf a}=k_1+_i(k_2\cdot_j{\bf a});$

\no provided all those operation results exist, then
$\widetilde{V}$ is called a multi-vector space on the multi-filed
$\widetilde{F}$ with a binary operation set $O(\widetilde{V})$,
denoted by $(\widetilde{V}; \widetilde{F})$.}

\vskip 3mm

Similar to multi-groups, results can be also obtained for
multi-rings or multi-vector spaces by generalizing classical
results in rings or linear spaces. Notice that in the references
$[17]$ and $[18]$, some such results have been gotten.

\vskip 4mm

\no{\bf $3.2$. The combinatorization of geometries}

\vskip 3mm

\no First, we generalize classical metric spaces by the
combinatorial speculation.

\vskip 4mm

\no{\bf Definition $3.5$} \ {\it A multi-metric space is a union
$\widetilde{M}=\bigcup\limits_{i=1}^m M_i$ such that each $M_i$ is
a space with metric $\rho_i$ for $\forall i, 1\leq i\leq m$.}

\vskip 3mm

Two well-known results in metric spaces are generalized.

\vskip 4mm

\no{\bf Theorem $3.3$}([19]) \ {\it Let
$\widetilde{M}=\bigcup\limits_{i=1}^m M_i$ be a completed
multi-metric space. For an $\epsilon$-disk sequence
$\{B(\epsilon_n,x_n)\}$, where $\epsilon_n > 0$ for $n=1,2,3,
\cdots$, the following conditions hold:}

($i$) \ {\it $B(\epsilon_1,x_1)\supset B(\epsilon_2,x_2)\supset
B(\epsilon_3,x_3)\supset\cdots\supset
B(\epsilon_n,x_n)\supset\cdots$};

($ii$) \ {\it $\lim\limits_{n\to+\infty}\epsilon_n=0$}.

\no {\it Then $\bigcap\limits_{n=1}^{+\infty}B(\epsilon_n,x_n)$
only has one point. }

\vskip 4mm

\no{\bf Theorem $3.4$}([19]) \ {\it Let
$\widetilde{M}=\bigcup\limits_{i=1}^m M_i$ be a completed
multi-metric space and $T$ a contraction on $\widetilde{M}$. Then}

$$1\leq ^{\#}\Phi (T)\leq m,$$

\no {\it where $^{\#}\Phi(T)$ denotes the number of fixed points
of $T$.}

\vskip 3mm

Particularly, let $m=1$. We get the {\it Banach fixed-point
theorem} again.

\vskip 4mm

\no{\bf Corollary $3.1$}(Banach) \ {\it Let $M$ be a metric space
and $T$ a contraction on $M$. Then $T$ has just one fixed point.}

\vskip 3mm

{\it Smarandache geometries} were proposed by Smarandache in
$[25]$ which are generalization of classical geometries, i.e.,
these {\it Euclid, Lobachevshy-Bolyai-Gauss} and {\it Riemann
geometries} may be united altogether in a same space, by some
Smarandache geometries under the combinatorial speculation. These
geometries can be either partially Euclidean and partially
Non-Euclidean, or Non-Euclidean. In general, Smarandache
geometries are defined by the next definition.

\vskip 4mm

\no{\bf Definition $3.6$} {\it An axiom is said to be
Smarandachely denied if the axiom behaves in at least two
different ways within the same space, i.e., validated and
invalided, or only invalided but in multiple distinct ways.

A Smarandache geometry is a geometry which has at least one
Smarandachely denied axiom($1969$).}

\vskip 3mm

For example, let us consider an euclidean plane ${\bf R}^2$ and
three non-collinear points $A,B$ and $C$. Define $s$-points as all
usual euclidean points on ${\bf R}^2$ and $s$-lines as any
euclidean line that passes through one and only one of points
$A,B$ and $C$. Then this geometry is a Smarandache geometry
because two axioms are Smarandachely denied comparing with an
Euclid geometry:

($i$) \ The axiom (A5) that through a point exterior to a given
line there is only one parallel passing through it is now replaced
by two statements: {\it one parallel} and {\it no parallel}. Let
$L$ be an $s$-line passing through $C$ and is parallel in the
euclidean sense to $AB$. Notice that through any $s$-point not
lying on $AB$ there is one $s$-line parallel to $L$ and through
any other $s$-point lying on $AB$ there are no $s$-lines parallel
to $L$ such as those shown in Fig.$3.1(a)$.

\includegraphics[bb=10 10 200 120]{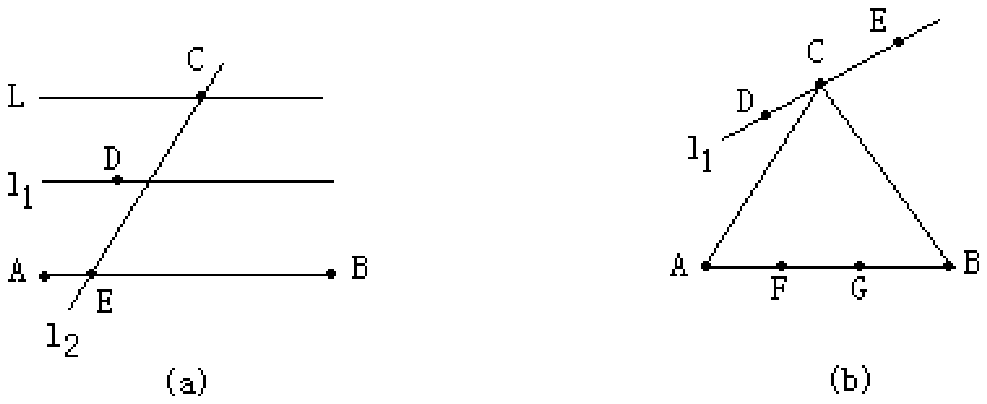}

\vskip 2mm

\c{\bf Fig.$3.1$} \vskip 2mm

($ii$) \ The axiom that through any two distinct points there
exists one line passing through them is now replaced by; {\it one
$s$-line} and {\it no $s$-line}. Notice that through any two
distinct $s$-points $D,E$ collinear with one of $A,B$ and $C$,
there is one $s$-line passing through them and through any two
distinct $s$-points $F,G$ lying on $AB$ or non-collinear with one
of $A,B$ and $C$, there is no $s$-line passing through them such
as those shown in Fig.$3.1(b)$.\vskip 2mm

A {\it Smarandache $n$-manifold} is an $n$-dimensional manifold
that support a Smarandache geometry. Now there are many approaches
for constructing Smarandache manifolds in the case of $n=2$. A
general way is by the so called {\it map geometries} without or
with boundary underlying orientable or non-orientable maps
proposed in references $[14]$ and $[15]$ firstly.

\vskip 4mm

\no{\bf Definition $3.7$} \ {\it For a combinatorial map $M$ with
each vertex valency$\geq 3$, endow each vertex $u, u\in V(M)$ a
real number $\mu (u), 0 \ < \mu (u) \ < \ \frac{4\pi}{\rho_M(u)}$.
Call $(M,\mu)$ a map geometry without boundary, $\mu:
V(M)\rightarrow R$ an angle function on $M$.}

\vskip 4mm

\no{\bf Definition $3.8$} \ {\it For a map geometry $(M,\mu )$
without boundary and faces $f_1,f_2,\cdots ,f_l$ $\in F(M), 1\leq
l\leq \phi (M)-1$, if $S(M)\setminus\{f_1,f_2,\cdots ,f_l\}$ is
connected, then call $(M,\mu)^{-l}= (S(M)\setminus\{f_1,f_2,\cdots
,f_l\}, \mu)$ a map geometry with boundary $f_1,f_2,\cdots ,f_l$,
where $S(M)$ denotes the locally orientable surface underlying
$M$.}

\vskip 3mm

A realization for vertices $u,v,w\in V(M)$ in a space ${\bf R}^3$
is shown in Fig.$3.2$, where $\rho_M(u)\mu (u) <  2\pi$,
$\rho_M(v)\mu (v) = 2\pi$ and $\rho_M(w)\mu (w)
> 2\pi$, are called to be elliptic, euclidean or
hyperbolic, respectively.

\vskip 2mm

\includegraphics[bb=30 20 100 110]{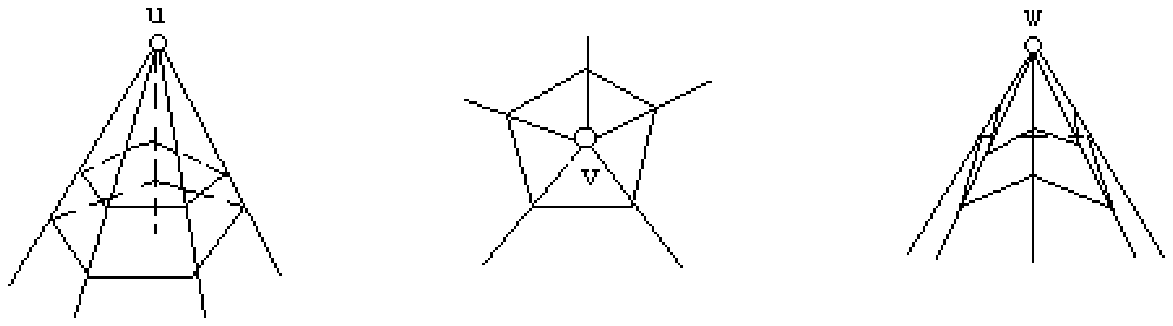}

\vskip 2mm

\c{$\rho_M(u)\mu (u) \ < \ 2\pi$\hskip 18mm $\rho_M(u)\mu (u)
=2\pi$\hskip 20mm $\rho_M(u)\mu (u) \ > \ 2\pi$}\vskip 4mm

\c{\bf Fig.$3.2$}\vskip 2mm

On an Euclid plane ${\bf R}^2$, a straight line passing through an
elliptic or a hyperbolic point is shown in Fig.$3.3$.

\includegraphics[bb=40 5 100 110]{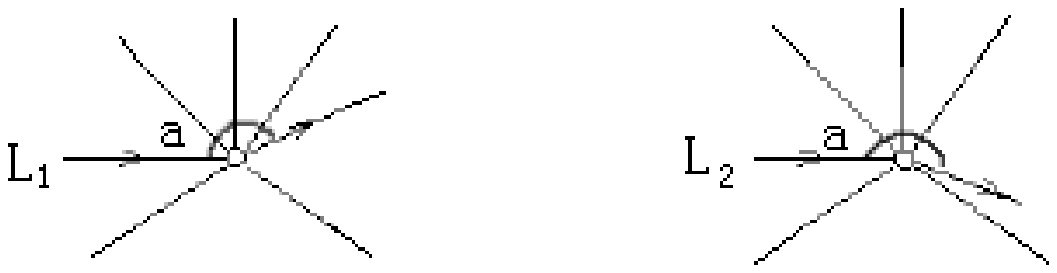}

\vskip 2mm

\c{\bf Fig.$3.3$}\vskip 4mm

\no{\bf Theorem $3.5$}([17])\ {\it There are Smarandache
geometries, including paradoxist geometries, non-geometries and
anti-geometries in map geometries without or with boundary.}

\vskip 3mm

Generally, we can generalize the ideas in Definitions $3.7$ and
$3.8$ to metric spaces further and find new geometries.

\vskip 4mm

\no{\bf Definition $3.9$} \ {\it Let $U$ and $W$ be two metric
spaces with metric $\rho$, $W\subseteq U$. For $\forall u\in U$,
if there is a continuous mapping $\omega: u\rightarrow \omega(u)$,
where $\omega(u)\in {\bf R}^n$ for an integer $n, n\geq 1$ such
that for any number $\epsilon >0$, there exist a number $\delta
>0$ and a point $v\in W$, $\rho(u-v)<\delta$ such that
$\rho(\omega(u)-\omega(v))<\epsilon$, then $U$ is called a metric
pseudo-space if $U=W$ or a bounded metric pseudo-space if there is
a number $N > 0$ such that $\forall w\in W$, $\rho(w)\leq N$,
denoted by $(U,\omega)$ or $(U^-,\omega)$, respectively.}

\vskip 3mm

For the case $n=1$, we can also explain $\omega(u)$ being an angle
function with $0 <\omega(u)\leq 4\pi$ as the case in map
geometries without or with boundary, i.e.,

\[
\omega(u)=\left\{
\begin{array}{ll}
\omega(u)(mod 4\pi), & {\rm if } \  {\rm u\in W},\\
2\pi, & {\rm if } \  {\rm u\in U\setminus W}
\end{array}
\right.
\]

\no and get some interesting metric pseudo-space geometries. For
example, let $U=W={\rm Euclid \ plane}=\sum$, then we obtained
some interesting results for pseudo-plane geometries $(\sum,
\omega)$ as shown in the following([17]).

\vskip 4mm

\no{\bf Theorem $3.6$} \ {\it In a pseudo-plane $(\sum,\omega)$,
if there are no euclidean points, then all points of
$(\sum,\omega)$ is either elliptic or hyperbolic.}

\vskip 4mm

\no{\bf Theorem $3.7$} \ {\it There are no saddle points and
stable knots in a pseudo-plane $(\sum,\omega)$.}

\vskip 4mm

\no{\bf Theorem $3.8$} \ {\it For two constants $\rho_0,\theta_0$,
$\rho_0 >0$ and $\theta_0\not= 0$, there is a pseudo-plane
$(\sum,\omega)$ with}

$$\omega(\rho,\theta)=2( \pi - \frac{\rho_0}{\theta_0\rho})\ {or} \
 \omega(\rho,\theta)=2(\pi+ \frac{\rho_0}{\theta_0\rho})$$

\no{\it such that}

$$\rho = \rho_0$$

\no{\it is a limiting ring in $(\sum,\omega)$.}

\vskip 3mm

Now for an $m$-manifold $M^m$ and $\forall u\in M^m$, choose
$U=W=M^m$ in Definition $3.9$ for $n=1$ and $\omega(u)$ a smooth
function. We get pseudo-manifold geometries $(M^m,\omega)$ on
$M^m$. By the reference $[2]$, a {\it Minkowski norm} on $M^m$ is
a function $F:M^m\rightarrow [0,+\infty)$ such that\vskip 3mm

($i$) \ \ \ $F$ is smooth on $M^m\setminus\{0\}$;

($ii$) \ \ $F$ is $1$-homogeneous, i.e.,
$F(\lambda\overline{u})=\lambda F(\overline{u})$ for
$\overline{u}\in M^m$ and $\lambda >0$;

($iii$) \ for $\forall y\in M^m\setminus\{0\}$, the symmetric
bilinear form $g_y: M^m\times M^m\rightarrow R$ with

$$g_y(\overline{u},\overline{v})=\frac{1}{2}
\frac{\partial^2F^2(y+s\overline{u}+t\overline{v})}{\partial
s\partial t}|_{t=s=0}$$\vskip 2mm

\no is positive definite and a {\it Finsler manifold} is a
manifold $M^m$ endowed with a function $F:TM^m\rightarrow
[0,+\infty)$ such that\vskip 3mm

($i$) \ \ $F$ is smooth on
$TM^m\setminus\{0\}=\bigcup\{T_{\overline{x}}M^m\setminus\{0\}:\overline{x}\in
M^m \}$;

($ii$) \ $F|_{T_{\overline{x}}M^m}\rightarrow[0,+\infty)$ is a
Minkowski norm for $\forall\overline{x}\in M^m$.\vskip 2mm

As a special case, we choose $\omega(\overline{x})=
F(\overline{x})$ for $\overline{x}\in M^m$, then $(M^m,\omega)$ is
a {\it Finsler manifold}. Particularly, if
$\omega(\overline{x})=g_{\overline{x}}(y,y)= F^2(x,y)$, then
$(M^m,\omega)$ is a {\it Riemann manifold}. Therefore, we get a
relation for Smarandache geometries with Finsler or Riemann
geometry.

\vskip 4mm

\no{\bf Theorem $3.9$} \ {\it There is an inclusion for
Smarandache, pseudo-manifold, Finsler and Riemann geometries as
shown in the following:}

\begin{eqnarray*}
\{Smarandache \ geometries\}&\supset&
\{pseudo-manifold \ geometries\}\\
&\supset&\{Finsler \ geometry\}\\
&\supset&\{Riemann \ geometry\}.
\end{eqnarray*}

\vskip 5mm

\no{\bf $4$. The contribution of combinatorial speculation to
theoretical physics}

\vskip 4mm

\no The progress of theoretical physics in last twenty years of
the 20th century enables human beings to probe the mystic cosmos:
{\it where are we came from? where are we going to?} Today, these
problems still confuse eyes of human beings. Accompanying with
research in cosmos, new puzzling problems also arose: {\it Whether
are there finite or infinite cosmoses? Is just one? What is the
dimension of our cosmos?} {\it We do not even know what the right
degree of freedom in the universe is}, as Witten said([3]).

We are used to the idea that our living space has three
dimensions: {\it length, breadth} and {\it height}, with time
providing the fourth dimension of spacetime by Einstein. Applying
his {\it principle of general relativity}, i.e. {\it all the laws
of physics take the same form in any reference system} and the
{\it equivalence principle}, i.e., {\it there are no difference
for physical effects of the inertial force and the gravitation in
a field small enough}, Einstein got the {\it equation of
gravitational field}

$$R_{\mu\nu}-\frac{1}{2}Rg_{\mu\nu}+\lambda g_{\mu\nu}=-8\pi GT_{\mu\nu}.$$

\no where $R_{\mu\nu}=R_{\nu\mu}=R_{\mu i\nu}^{\alpha}$,

$$R_{\mu i\nu}^{\alpha}=\frac{\partial\Gamma_{\mu i}^i}{\partial x^{\nu}}-
\frac{\partial\Gamma_{\mu\nu}^i}{\partial x^{i}}+\Gamma_{\mu
i}^{\alpha}\Gamma_{\alpha\nu}^i-\Gamma_{\mu\nu}^{\alpha}\Gamma_{\alpha
i}^i,$$

$$\Gamma_{mn}^g=\frac{1}{2}g^{pq}(\frac{\partial g_{mp}}{\partial u^n}
+\frac{\partial g_{np}}{\partial u^m}-\frac{\partial
g_{mn}}{\partial u^p})$$

\no and $R=g^{\nu\mu}R_{\nu\mu}$.

Combining the Einstein's equation of gravitational field with the
{\it cosmological principle}, i.e., {\it there are no difference
at different points and different orientations at a point of a
cosmos on the metric $10^4l.y.$} , Friedmann got a standard model
of cosmos. The metrics of the standard cosmos are

$$ds^2 = -c^2dt^2+a^2(t)[\frac{dr^2}{1-Kr^2}+r^2(d\theta^2+\sin^2\theta d\varphi^2)]$$

\no and

$$g_{tt}=1, \ g_{rr}=-\frac{R^2(t)}{1-Kr^2}, g_{\phi\phi}=-r^2R^2(t)\sin^2\theta.$$

The standard model of cosmos enables the birth of big bang model
for our cosmos in thirties of the 20th century. The following
diagram describes the developing process of our cosmos in
different periods after the big bang.

\includegraphics[bb=-40 5 100 270]{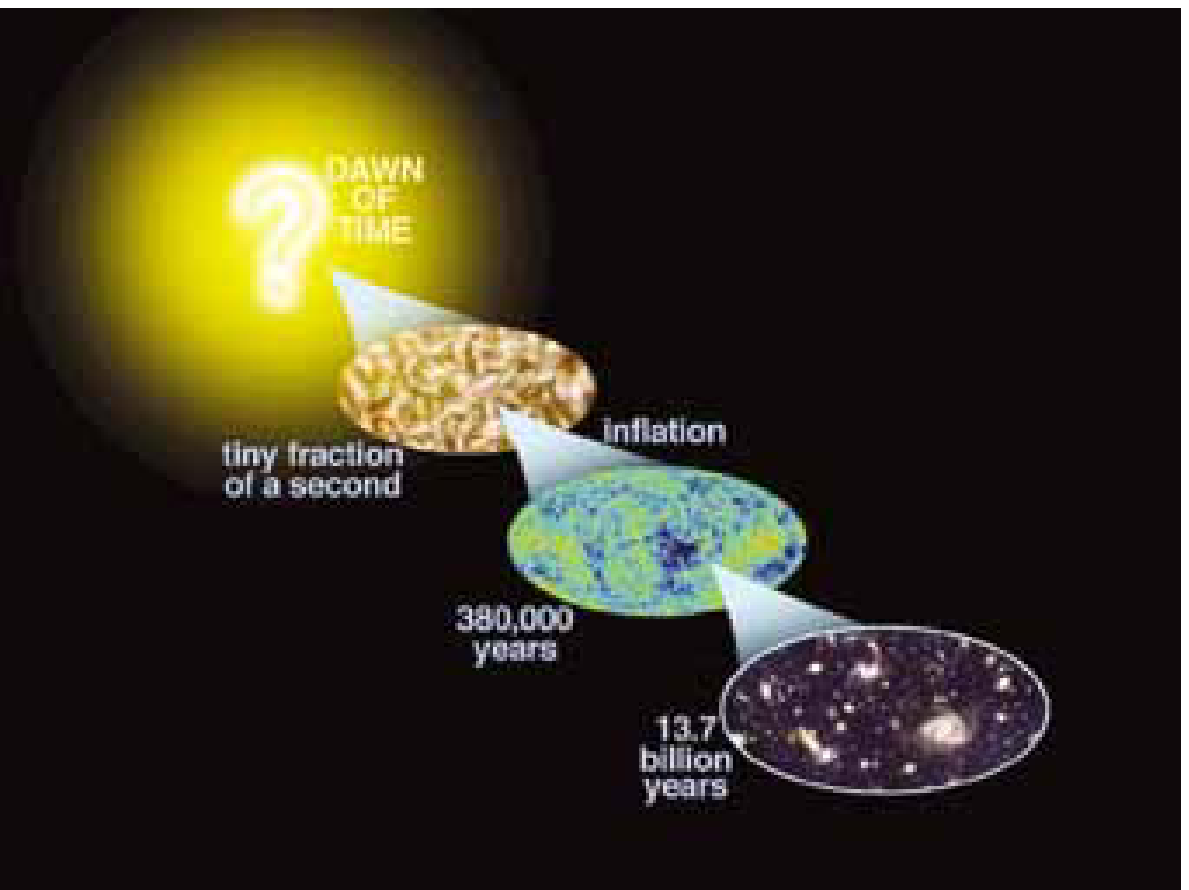}

\vskip 3mm

\c{\bf Fig.$4.1$}

\vskip 4mm

\no{\bf $4.1$. The M-theory}

\vskip 3mm

\no The M-theory was established by Witten in 1995 for the unity
of those five already known string theories and superstring
theories, which postulates that all matter and energy can be
reduced to {\it branes} of energy vibrating in an $11$ dimensional
space, and then in a higher dimensional space solve the Einstein's
equation of gravitational field under some physical conditions
([1],[3],[22]-[23]). Here, a {\it brane} is an object or subspace
which can have various spatial dimensions. For any integer $p\geq
0$, a {\it $p$-brane} has length in $p$ dimensions. For example, a
{\it $0$-brane} is just a point or particle; a {\it $1$-brane} is
a string and a {\it $2$-brane} is a surface or membrane, $\cdots$.

One mainly discuss line elements in differential forms in Riemann
geometry. By a geometrical view, these $p$-branes in M-theory can
be seen as these volume elements in spaces in the counterpart.
Whence, we can construct a graph model for $p$-branes in a space
and combinatorially research graphs in spaces.

\vskip 4mm

\no{\bf Definition $4.1$} \ {\it For each $m$-brane ${\bf B}$ of a
space ${\bf R}^{m}$, let $(n_1({\bf B}),n_2({\bf B}),\cdots,
n_p({\bf B}))$ be its unit vibrating normal vector along these $p$
directions and $q:{\bf R}^{m}\rightarrow{\bf R}^4$ a continuous
mapping. Now construct a graph phase $({\mathcal
G},\omega,\Lambda)$ by}

$$V({\mathcal G})=\{p-branes \ q({\bf B})\},$$

$$E({\mathcal G})=\{(q({\bf B}_1),q({\bf B}_2))| there \ is \ an \ action \ between \ {\bf B}_1 \
and \ {\bf B}_2\},$$

$$\omega(q({\bf B}))=(n_1({\bf B}),n_2({\bf B}),\cdots,
n_p({\bf B})),$$

\no {\it and}

$$\Lambda(q({\bf B}_1),q({\bf B}_2))= \ forces \ between \ {\bf B}_1 \ and \ {\bf B}_2.$$

\no{\it Then we get a graph phase $({\mathcal G},\omega,\Lambda)$
in ${\bf R}^4$. Similarly, if $m=11$, it is a graph phase for the
M-theory.}

\vskip 3mm

As an example for applying M-theory to find an accelerating
expansion cosmos of $4$-dimensional cosmos from supergravity
compactification on hyperbolic spaces is the {\it
Townsend-Wohlfarth type metric} in which the line element is

$$ds^2=e^{-m\phi(t)}(-S^6dt^2+S^2dx_3^2)+r_C^2e^{2\phi(t)}ds_{H_m}^2,$$

\no where

$$\phi(t)=\frac{1}{m-1}(\ln K(t)-3\lambda_0t), \ \ \
S^2=K^{\frac{m}{m-1}}e^{-\frac{m+2}{m-1}\lambda_0t}$$

\no and

$$K(t)=\frac{\lambda_0\zeta r_c}{(m-1)\sin[\lambda_0\zeta|t+t_1|]}$$

\no with $\zeta=\sqrt{3+6/m}$. This solution is obtainable from
space-like brane solution and if the proper time $\varsigma$ is
defined by $d\varsigma= S^3(t)dt$, then the conditions for
expansion and acceleration are $\frac{dS}{d\varsigma}>0$ and
$\frac{d^2S}{d\varsigma^2}>0$. For example, the expansion factor
is $3.04$ if $m=7$, i.e., a really expanding cosmos.

According to M-theory, the evolution picture of our cosmos started
as a perfect $11$ dimensional space. However, this $11$
dimensional space was unstable. The original $11$ dimensional
space finally cracked into two pieces, a $4$ and a $7$ dimensional
cosmos. The cosmos made the $7$ of the $11$ dimensions curled into
a tiny ball, allowing the remaining $4$ dimensional cosmos to
inflate at enormous rates.

\vskip 4mm

\no{\bf $4.2$. The combinatorial cosmos}

\vskip 3mm

\no The combinatorial speculation enables us to introduce the
conception of combinatorial cosmoss([17]).

\vskip 4mm

\no{\bf Definition $4.2$} \ {\it A combinatorial cosmos is
constructed by a triple $(\Omega,\Delta,T)$, where }

$$\Omega=\bigcup\limits_{i\geq 0}\Omega_i, \  \ \ \Delta=\bigcup\limits_{i\geq 0}O_i$$

\no{\it and $T=\{t_i; i\geq 0\}$ are respectively called the
cosmos, the operation or the time set with the following
conditions hold.}\vskip 2mm

($1$) \ {\it $(\Omega,\Delta)$ is a Smarandache multi-space
dependent on $T$, i.e., the cosmos $(\Omega_i,O_i)$ is dependent
on time parameters $t_i$ for any integer $i, i\geq 0$.}

($2$) \ {\it For any integer $i, i\geq 0$, there is a sub-cosmos
sequence}

$$(S): \ \Omega_i\supset\cdots\supset\Omega_{i1}\supset\Omega_{i0}$$

\no{\it in the cosmos $(\Omega_i,O_i)$ and for two sub-cosmoses
$(\Omega_{ij},O_i)$ and $(\Omega_{il}, O_i)$, if
$\Omega_{ij}\supset\Omega_{il}$, then there is a homomorphism
$\rho_{\Omega_{ij},\Omega_{il}}:(\Omega_{ij},O_i)\rightarrow(\Omega_{il},O_i)$
such that}\vskip 2mm

($i$) \ {\it for $\forall (\Omega_{i1},O_i),(\Omega_{i2},O_i),
(\Omega_{i3},O_i)\in(S)$, if
$\Omega_{i1}\supset\Omega_{i2}\supset\Omega_{i3}$, then}

$$\rho_{\Omega_{i1},\Omega_{i3}}=
\rho_{\Omega_{i1},\Omega_{i2}}\circ\rho_{\Omega_{i2},\Omega_{i3}},$$

\no{\it where ¡°$\circ$¡± denotes the composition operation on
homomorphisms.}

($ii$) \ {\it for $\forall g,h\in\Omega_i$, if for any integer
$i$, $\rho_{\Omega,\Omega_i}(g)=\rho_{\Omega,\Omega_i}(h)$, then
$g=h$.}

($iii$) \ {\it for $\forall i$, if there is an $f_i\in\Omega_i$
with}

$$\rho_{\Omega_i,\Omega_i\bigcap\Omega_j}(f_i)=
\rho_{\Omega_j,\Omega_i\bigcap\Omega_j}(f_j)$$

\no{\it for integers $i,j,\Omega_i\bigcap\Omega_j\not=\emptyset$,
then there exists an $f\in\Omega$ such that
$\rho_{\Omega,\Omega_i}(f)=f_i$ for any integer $i$.}

\vskip 3mm

By this definition, there is just one cosmos $\Omega$ and the
sub-cosmos sequence is

$${\bf R}^4\supset{\bf R}^3\supset
{\bf R}^2\supset{\bf R}^1\supset{\bf R}^0=\{P\}\supset{\bf
R}_7^-\supset\cdots\supset{\bf R}_1^-\supset{\bf R}_0^-=\{Q\}.$$

\no in the string/M-theory. In Fig.$4.1$, we show a combinatorial
cosmos with a higher dimensional cosmos outside our visual cosmos.

\includegraphics[bb=-60 5 100 180]{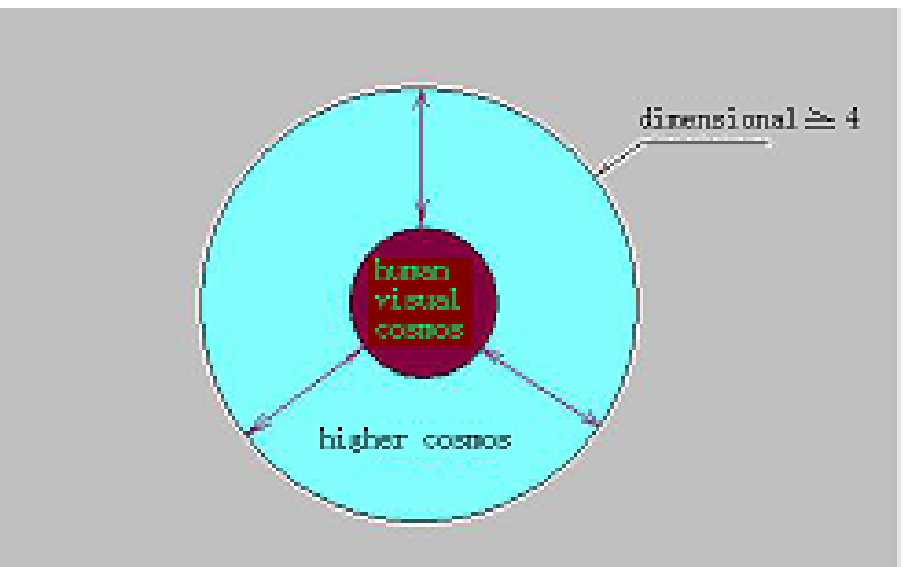}

\vskip 3mm

\c{\bf Fig.$4.2$}\vskip 3mm

If the dimensional of this cosmos outside our visual cosmos is $5$
or $6$, there has been established a dynamical theory by this
combinatorial speculation([20][21]). In this dynamics, we look for
a solution in the Einstein equation of gravitational field in
$6$-dimensional spacetime with a metric of the form

$$ds^2=-n^2(t,y,z)dt^2+a^2(t,y,z)d\sum_k^2+b^2(t,y,z)dy^2+d^2(t,y,z)dz^2$$

\no where $d\sum_k^2$ represents the $3$-dimensional spatial
sections metric with $k=-1,0,1$ respective corresponding to the
hyperbolic, flat and elliptic spaces. For $5$-dimensional
spacetime, deletes the undefinite $z$ in this metric form. Now
consider a $4$-brane moving in a $6$-dimensional {\it
Schwarzschild-ADS spacetime}, the metric can be written as

$$ds^2=-h(z)dt^2+\frac{z^2}{l^2}d\sum_k^2+h^{-1}(z)dz^2,$$

\no where

$$d\sum_k^2=\frac{dr^2}{1-kr^2}+r^2d\Omega^2_{(2)}+(1-kr^2)dy^2$$

\no and

$$h(z)=k+\frac{z^2}{l^2}-\frac{M}{z^3}.$$

Then the equation of a $4$-dimensional cosmos moving in a
$6$-spacetime is

$$2\frac{\ddot{R}}{R}+3(\frac{\dot{R}}{R})^2=-3\frac{\kappa_{(6)}^4}{64}\rho^2-
\frac{\kappa_{(6)}^4}{8}\rho p-3\frac{\kappa}{R^2}-\frac{5}{l^2}$$

\no by applying the {\it Darmois-Israel conditions} for a moving
brane. Similarly, for the case of $a(z)\not=b(z)$, the equations
of motion of the brane are

$$\frac{d^2\dot{d}\dot{R}-d\ddot{R}}{\sqrt{1+d^2\dot{R}^2}}-
\frac{\sqrt{1+d^2\dot{R}^2}}{n}(d\dot{n}\dot{R}
+\frac{\partial_zn}{d}-(d\partial_zn-n\partial_zd)\dot{R}^2)
=-\frac{\kappa_{(6)}^4}{8}(3(p+\rho)+\hat{p}),$$

$$\frac{\partial_za}{ad}\sqrt{1+d^2\dot{R}^2}=-\frac{\kappa_{(6)}^4}{8}(\rho+p-\hat{p}),$$

$$\frac{\partial_zb}{bd}\sqrt{1+d^2\dot{R}^2}=-\frac{\kappa_{(6)}^4}{8}(\rho-3(p-\hat{p})),$$

\no where the energy-momentum tensor on the brane is

$$\hat{T}_{\mu\nu}=h_{\nu\alpha}T_{\mu}^{\alpha}-\frac{1}{4}Th_{\mu\nu}$$

\no with $T_{\mu}^{\alpha}= diag(-\rho,p,p,p,\hat{p})$ and the
{\it Darmois-Israel conditions}

$$[K_{\mu\nu}]=-\kappa_{(6)}^2\hat{T}_{\mu\nu},$$

\no where $K_{\mu\nu}$ is the extrinsic curvature tensor.

The idea of combinatorial cosmoses also presents new questions to
combinatorics, such as:\vskip 3mm

($i$)\ \ \ to embed a graph into spaces with dimensional$\geq 4$;

($ii$)\ \ to research the phase space of a graph embedded in a
space;

($iii$)\ to establish graph dynamics in a space  with
dimensional$\geq 4$, $\cdots$, etc..\vskip 2mm

For example, we have gotten the following result for graphs in
spaces in [17].

\vskip 4mm

\no{\bf Theorem $4.1$} \ {\it A graph $G$ has a nontrivial
including multi-embedding on spheres $P_1\supset
P_2\supset\cdots\supset P_s$ if and only if there is a block
decomposition $G=\biguplus\limits_{i=1}^sG_i$ of $G$ such that for
any integer $i, 1 < i < s$,}

($i$) \ {\it $G_i$ is planar;}

($ii$) \ {\it for $\forall v\in V(G_i)$, $N_G(x)\subseteq
(\bigcup\limits_{j=i-1}^{i+1}V(G_j))$.}

\vskip 3mm

Further research of these combinatorial cosmoses will richen the
knowledge of combinatorics and cosmology, also get the
combinatorization for cosmology.

\newpage

\no{\bf References}

\vskip 6mm

\re{[1]}I.Antoniadis, Physics with large extra dimensions: String
theory under experimental test, {\it Current Science}, Vol.81,
No.12, 25(2001),1609-1613

\re{[2]}S.S.Chern and W.H.Chern, {\it Lectures in Differential
Geometry}(in Chinese), Peking University Press, 2001.

\re{[3]}M.J.Duff, A layman's guide to M-theory, {\it arXiv}:
hep-th/9805177, v3, 2 July(1998).

\re{[4]}B.J.Fei, {\it Relativity Theory and Non-Euclid
Geometries}, Science Publisher Press, Beijing, 2005.

\re{[5]}U.G$\ddot{u}$nther and A.Zhuk, Phenomenology of
brane-world cosmological models, {\it arXiv: gr-qc/0410130}.

\re{[6]}S.Hawking, {\it A Brief History of Times}, A Bantam
Books/November, 1996.

\re{[7]}S.Hawking, {\it The Universe in Nutshell}, A Bantam
Books/November, 2001.

\re{[8]}D.Ida, Brane-world cosmology, {\it arXiv: gr-qc/9912002}.

\re{[9]}H.Iseri, {\it Smarandache Manifolds}, American Research
Press, Rehoboth, NM,2002.

\re{[10]}H.Iseri, {\it Partially Paradoxist Smarandache
Geometries}, http://www.gallup.unm.
edu/\~smarandache/Howard-Iseri-paper.htm.

\re{[11]}M.Kaku, {\it Hyperspace: A Scientific Odyssey through
Parallel Universe, Time Warps and 10th Dimension}, Oxford Univ.
Press.

\re{[12]}P.Kanti, R.Madden and K.A.Olive, A 6-D brane world model,
{\it arXiv: hep-th/0104177}.

\re{[13]}L.Kuciuk and M.Antholy, An Introduction to Smarandache
Geometries, {\it Mathematics Magazine, Aurora, Canada},
Vol.12(2003)

\re{[14]}L.F.Mao, {\it On Automorphisms groups of Maps, Surfaces
and Smarandache geometries}, {\it Sientia Magna}, Vol.$1$(2005),
No.$2$, 55-73.

\re{[15]}L.F.Mao, A new view of combinatorial maps by
Smarandache's notion, {\it arXiv: math.GM/0506232}.

\re{[16]}L.F.Mao, {\it Automorphism Groups of Maps, Surfaces and
Smarandache Geometries}, American Research Press, 2005.

\re{[17]}L.F.Mao, {\it Smarandache multi-space theory}, Hexis,
Phoenix, AZ£¬2006.

\re{[18]}L.F.Mao, On algebraic multi-groups, {\it Sientia Magna},
Vol.$2$(2006), No.$1$, 64-70.

\re{[19]}L.F.Mao, On multi-metric spaces, {\it Sientia Magna},
Vol.$2$(2006), No.$1$, 87-94.

\re{[20]}E.Papantonopoulos, Braneworld cosmological models, {\it
arXiv: gr-qc/0410032}.

\re{[21]}E.Papantonopoulos, Cosmology in six dimensions, {\it
arXiv: gr-qc/0601011}.

\re{[22]}J.A.Peacock, {\it Cosmological Physics}, Cambridge
University Press, 2003.

\re{[23]}J.Polchinski, {\it String Theory}, Vol.1-Vol.2, Cambridge
University Press, 2003.

\re{[24]}D.Rabounski, Declaration of academic freedom: Scientific
Human Rights, {\it Progress in Physics}, January (2006), Vol.1,
57-60.

\re{[25]}F.Smarandache, Mixed noneuclidean geometries, {\it eprint
arXiv: math/0010119}, 10/2000.

\re{[26]}F.Smarandache, {\it A Unifying Field in Logics.
Neutrosopy: Neturosophic Probability, Set, and Logic}, American
research Press, Rehoboth, 1999.

\re{[27]}F.Smarandache, Neutrosophy, a new Branch of Philosophy,
{\it Multi-Valued Logic}, Vol.8, No.3(2002)(special issue on
Neutrosophy and Neutrosophic Logic), 297-384.

\re{[28]}F.Smarandache, A Unifying Field in Logic: Neutrosophic
Field, {\it Multi-Valued Logic}, Vol.8, No.3(2002)(special issue
on Neutrosophy and Neutrosophic Logic), 385-438.

\re{[29]}J.Stillwell, {\it Classical topology and combinatorial
group theory}, Springer-Verlag New York Inc., (1980).

\end{document}